\documentclass[11pt]{amsart}
\usepackage{amsfonts,amssymb,epsf}

\def\Pr{\begin{proof}}
\def\Rp{\end{proof}}

\def\ln{\left\Vert}
\def\rn{\right\Vert}
\newcommand\norme[1]{\|#1\|}
\def\IR{\mathbb R}
\def\IN{\mathbb N}

\def\ZF{\mathbf {ZF}}
\def\DC{\mathbf {DC}}
\def\HB{\mathbf {HB}}

\def\AC{\mathbf {AC}}

\def\ZFC{\mathbf {ZFC}}

\theoremstyle{plain}

\newtheorem{corollary}{Corollary}

\newtheorem*{proposition*}{Proposition}
\newtheorem{theorem}{Theorem}
\newtheorem*{theorem*}{Theorem}
\newtheorem*{corollary*}{Corollary}
\newtheorem{lemma}{Lemma}

\theoremstyle{definition}

\newtheorem{notation}{Notation}
\newtheorem{question}{Question}

\theoremstyle{remark}
\newtheorem{remark}{Remark}

\newtheorem*{axiom}{}

\begin{document}
\title
{A new proof of James' sup theorem}
\author[M.~Morillon]{Marianne Morillon}
\address{ERMIT, D\'epartement de Math\'ematiques et Informatique,
Universit\'e de La R\'eunion, 15 avenue Ren\'e Cassin - BP 7151 -
97715 Saint-Denis Messag. Cedex 9 FRANCE}
\email[Marianne Morillon]{mar@univ-reunion.fr}
\urladdr{http://univ-reunion.fr/$\sim$mar}

\begin{abstract} We provide  a new proof of James' sup theorem for (non necessarily separable)
Banach spaces. One of the ingredients is the following generalization of a theorem
 of Hagler and Johnson (\cite{Ha-Jo})~: {\em ``If a normed space $E$ does not contain any
asymptotically isometric copy of $\ell^1(\IN)$, then every bounded sequence of $E'$
has a normalized block sequence pointwise converging to $0$''}.
\end{abstract}

\maketitle

\bigskip

\begin{center}
   \mbox{\epsfbox{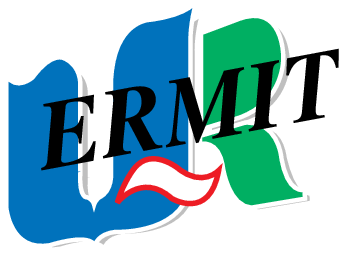}}
   \end{center}

\begin{center}
EQUIPE R\'EUNIONNAISE DE MATH\'EMATIQUES ET INFORMATIQUE TH\'EORIQUE
\end{center}
\bigskip

\bigskip

\begin{center}
http://univ-reunion.fr/$\sim$dirermit
\end{center}

\

\bigskip

\begin{center}
Preprint 7 janvier  2005
\end{center}

\newpage

\section{Preliminaries}
\subsection{Introduction}
Given a normed space $E$,  we denote  by $\Gamma_E$ its closed unit ball,
by $E'$ the continuous dual of $E$ (endowed with its dual norm),
 and by $E''$ the second continuous dual of $E$.
 Say  that  the normed space $E$ is {\em onto-reflexive} if the canonical mapping $j_E : E \to E''$
is onto.  Say that $E$  is {\em J-reflexive} (James-reflexive) if $\Gamma_E$
does not contain any
 sequence $(a_n)_{n \in \IN}$ satisfying
 $$\inf_{n \in \IN} distance(span\{a_i : \; i < n\}, conv\{a_i : i \ge n\}) >0$$
    Say that   $E$ is {\em sup-reflexive}
if every $f \in E'$ attains its upper bound on $\Gamma_E$.
It is known that, with the  Axiom of Choice,
 for a given Banach space the following notions are
equivalent~:  onto-reflexivity, weak compactness of the closed unit ball,
J-reflexivity, sup-reflexivity,
Eberlein-Smulyan property (``Every bounded sequence has a weakly convergent subsequence'')~\dots \;
Both implications ``J-reflexivity $\Rightarrow$ onto-reflexivity''
and ``sup-reflexivity $\Rightarrow$ J-reflexivity''
are due to James~; the first implication has a recent short  proof (see \cite{Oj98}), but
classical proofs of
the ``reflexivity'' of sup-reflexive Banach spaces  (James' sup theorem) are rather intricate
(see~\cite{Ja64b}, \cite{Ja72}, \cite{Pr}).
There exist  simpler proofs  of this theorem {\em under some restrictions}~: for example,
there is a short proof, relying on Simons' inequality, that
{\em separable} sup-reflexive Banach spaces are onto-reflexive (see~\cite{DGZ}).
In this paper, we provide a new proof of James' sup theorem.

\subsection{Presentation of the results}
We work in Zermelo-Fraenkel set-theory without choice $\ZF$,
and expressly mention the
two (weak) forms of the Axiom of Choice
($\AC$) that are used in our proofs, namely~:
the axiom of {\em Hahn-Banach} ($\HB$)
and the axiom of {\em Dependent Choices} ($\DC$) (see Section~\ref{sec:AC}).

Given a real vector space $E$ and a sequence $(x_n)_{n \in \IN}$ of $E$, 
a sequence $(b_n)_{n \in \IN}$ of
$E$ is a {\em block sequence} of $(x_n)_n$ if there exists a sequence $(F_n)_{n \in \IN}$ of
pairwise disjoint finite subsets of $\IN$
and a sequence $(\lambda_i)_{i \in \IN}$ of real numbers
such that for every $n \in \IN$, $b_n = \sum_{i \in F_n} \lambda_i x_i$~;
 if for each $n \in \IN$, $\sum_{i \in F_n}|\lambda_i|=1$, the block sequence $(b_n)_{n \in \IN}$
is said to be {\em normalized}~; if, in addition,
for every $i \in \cup_{n \in \IN}F_n$, $\lambda_i \ge 0$,
say that  the  block sequence is {\em convex}~: thus, every infinite subsequence 
is a convex block sequence.
The sequence $(F_n)_n$ is called a  sequence of {\em supports} of the block sequence $(b_n)_n$.
Say that a topological space is {\em sequentially compact}
if every sequence of this space has an infinite subsequence which converges.
Say that a  subset $C$ of a topological vector space $E$ is {\em block compact}
if every sequence of $C$ has a normalized block sequence $(b_n)_{n \in \IN}$ which converges 
in $E$
(hence the normalized block sequence $(\frac{b_{2n} - b_{2n+1}}{2})_{n \in \IN}$ 
converges to $0$). Say that $C$
is {\em convex block compact} if every sequence of $C$ has a convex  block sequence which 
converges in $E$. \par
Using Simons' inequality,
we begin by proving   in $\ZF \mathbf + \HB$
that {\em every sup-reflexive normed space with a *weak convex block compact dual ball 
is J-reflexive} (Corollary~\ref{cor:sep}-\ref{it:c-b-c} of Section~\ref{sec:dual-ball-si})~; 
in particular, this
yields a new proof (indeed in $\ZF$) of the J-reflexivity of sup-reflexive  
separable normed spaces (Corollary~\ref{cor:sep}-\ref{it:sep}).
Then, in Section~\ref{subsec:*weak-db}, using Simons' inequality and Rosenthal's $\ell^1$-theorem,
we prove  the following ``generalization''~:
\begin{theorem}[$\HB$] \label{lem:*wdual} Given a sup-reflexive normed space $E$,
if  $\Gamma_{E'}$ is *weak block compact,  then $E$ is J-reflexive.
\end{theorem}

We then  prove (Section~\ref{sec:ha-jo-zfc})~:
\begin{theorem} [$\DC$] \label{theo:ha-jo-zfc} If  a normed space $E$  does not contain
any asymptotically isometric copy of $\ell^1(\IN)$, then, $\Gamma_{E'}$ is *weak
block compact.
\end{theorem}
\noindent Here, say that $E$ contains an {\em asymptotically isometric copy of  $\ell^1(\IN)$}
if there  exists  a sequence $(a_n)_{n \in \IN}$
of $\Gamma_E$ and
  some sequence $(\delta_n)_{n \in \IN}$ of $]0,1[$ converging to $1$ satisfying
the following inequality for every finite sequence $(\lambda_i)_{0 \le i \le n}$ of $\IR$~:
$$\sum_{0 \le i \le n} \delta_i |\lambda_i|  \le \ln  \sum_{0 \le i \le n} \lambda_i a_i \rn$$
Theorem~\ref{theo:ha-jo-zfc} generalizes  a  result
due to Hagler and Johnson (\cite{Ha-Jo}, 1977),
where the normed space  contains no {\em isomorphic}
copy of $\ell^1(\IN)$.

 We finally prove (Section~\ref{sec:noasyiso})~:
\begin{theorem}[$\HB$] \label{theo:asy-l1} A sup-reflexive normed space
does not contain any asymptotically isometric copy of $\ell^1(\IN)$.
\end{theorem}
\noindent This result is a straightforward generalization
of a short theorem due to James (see~\cite[Theorem~2 p.~209]{Ja64c}).

 Thus we get the following  new proof  of James' sup theorem~:
\begin{corollary} [$\HB \mathbf + \DC$] \label{theo:jamessup} Every  sup-reflexive
normed space is J-reflexive.
\end{corollary}
\Pr Let $E$ be a sup-reflexive normed space. According to Theorem~\ref{theo:asy-l1},
the space  $E$ does not contain any asymptotically isometric copy of $\ell^1(\IN)$~; so, with
 Theorem~\ref{theo:ha-jo-zfc},   $\Gamma_{E'}$ is *weak block compact~; whence,
by Theorem~\ref{lem:*wdual},  $E$ is J-reflexive.
\Rp

\subsection{Questions}
Given a normed space $E$, obviously,
\begin{multline} \notag
\Gamma_{E'} \text{ *weak  sequentially compact} \Rightarrow \\
\Gamma_{E'} \text{ *weak convex block compact } \Rightarrow
\Gamma_{E'} \text{ *weak block compact}
\end{multline}
The first implication is not reversible in $\ZFC$ since
there exists a Banach space $E$ such that $\Gamma_{E'}$ endowed with the *weak topology is convex block compact and not
sequentially compact (see~\cite{Ha-Od})~: notice that the construction of the space built there 
depends on a well-order on $\IR$.
\begin{question} Does there exist in $\ZF$ a Banach space $E$ such that
$\Gamma_{E'}$ is *weak convex block compact and not *weak sequentially compact~?
\end{question}

\begin{question} Does there exist  (in $\ZF$ or $\ZFC$) a Banach space $E$ such that
$\Gamma_{E'}$ is *weak block compact and not *weak convex block compact~?
\end{question}

\begin{question} \label{ques:v-b-d-f-t} According to a theorem of  Valdivia
(see~\cite{Va}) for locally convex spaces,
who refers to Bourgain and Diestel (see \cite{Bou-Die}) for Banach spaces, themselves
referring to Bourgain-Fremlin-Talagrand (see~\cite{Bou-Fre-Tal}), $\DC$ implies that
{\em ``the dual ball of a  normed space not containing any isomorphic copy of $\ell^1(\IN)$
is *weak convex block compact''}. Does  this result
persist (in $\ZF$ or $\ZFC$) for normed spaces which
do not contain asymptotically isometric copies of
$\ell^1(\IN)$~?  
\end{question}

\medskip

\section{Spaces with a *weak block compact dual ball} \label{sec:dual-ball-si-ro}
\subsection{Spaces with a *weak convex block compact dual ball} \label{sec:dual-ball-si}
\begin{theorem*}[{\bf Simons' inequality} \cite{Si72}] Let $S$ be a set and $(x_n)_{n \in \IN}$
be a bounded sequence
of $\ell^{\infty}(S)$. Denote by $\Lambda$ the set of sequences
$(\lambda_n)_{n \in \IN} \in [0,1]^{\IN}$ satisfying  $\sum_{n \in \IN}\lambda_n=1$.
Assume that for every   $(\lambda_n)_{n \in \IN} \in \Lambda$,
the infinite convex combination $\sum_{n \in \IN}\lambda_n x_n$ attains its upper bound on $S$.
Then,
$$\inf \{ \sup_S \sum_{n \in \IN}\lambda_n x_n \; : \;
(\lambda_n)_{n \in \IN} \in \Lambda  \}
\le \sup_S \limsup_{n \in \IN} x_n  $$
\end{theorem*}
\noindent Notice that the proof given in \cite{Si72}  (see also \cite{oj98b}) is choiceless~: 
use convex 
combinations with  finite supports and rational coefficients.

\medskip

Given a normed space $E$, and some real number $\vartheta>0$,
say that a sequence $(a_n)_{n \in \IN}$ of $E$  is
a {\em $\vartheta$-sequence}  if
 $\inf_{n \in \IN} d(span\{a_i : \; i < n\}, conv\{a_i : i \ge n\}) \ge \vartheta$.
Given a  $\vartheta$-sequence $(a_n)_{n \in \IN}$ of $\Gamma_E$, and
denoting by $V$ the vector space $span(\{a_n : n \in \IN\})$,
  there is a sequence
 $(f_n)_{n \in \IN}$ of $\Gamma_{V'}$ satisfying $f_n(a_i)=0$ if $i<n$
 and $f_n(a_i) \ge \vartheta$ if $n \le i$~:
 in this case, say that  $(a_n,f_n)_{n \in \IN}$ 
 is  a {\em $\vartheta$-triangular} sequence of $E$.  Using $\HB$ (or rather its ``multiple form'',
 see Section~\ref{sec:AC}),  there exists a sequence $(\tilde f_n)_{n \in \IN}$ of $\Gamma_{E'}$
 such that each $\tilde f_n$
 extends $f_n$~: in this case, say that $(a_n,f_n)_{n \in \IN}$ is an {\em extended}
  triangular sequence of $E$. Thus, a normed space $E$ is J-reflexive
 if and only if it has  no $\vartheta$-triangular sequence  for any $\vartheta >0$.

\medskip

\begin{lemma}  \label{lem:sep-james} Given a sup-reflexive space $E$,
some $\vartheta>0$, and some extended $\vartheta$-triangular
sequence $(a_n,f_n)_{n \in \IN}$ of $E$,  no convex block sequence of
$(f_n)_{n \in \IN}$  pointwise converges.
\end{lemma}
\Pr Seeking a contradiction, assume that some convex block sequence $(b_n)_{n \in \IN}$
  pointwise converges to some $f$. Without loss of generality, we may assume that the
  sequence of supports $(F_n)_{n \in \IN}$ satisfies $F_0 < F_1 < \dots < F_n < \dots$.
 Observe that for every $n \in \IN$, $f(a_n)=0$.
 Then denoting by $h_n$ the mapping $\frac{b_n - f}{2}$, and by  $d_n$  the last element of $F_n$,
 the sequence
$(a_{d_n}, h_n)_{n \in \IN}$  is  $\frac{\vartheta}{2}$-triangular.
Using Simons' inequality and the assumption of sup-reflexivity  of $E$, 
there exists some finite convex combination
$g:=\sum_{n \in F} \lambda_n h_n$ of $(h_n)_{n \in \IN}$ such that
$\sup_{\Gamma_E} g \le \frac{\vartheta}{4}$~;
but for any integer $N > \max F$,  $g(a_{N}) =\sum_{i \in F} \lambda_i h_i(a_{N})
\ge \frac{\vartheta}{2}$~:  contradiction~!
\Rp


\begin{corollary} \label{cor:sep} Let $E$ be a sup-reflexive normed space.
\begin{enumerate}
\item \label{it:c-b-c}  If $\Gamma_{E'}$ is  *weak convex block compact,
 then $E$ has no extended triangular sequence, in which case, with $\HB$,
 $E$ is  J-reflexive.
\item \label{it:sep} If $E$ is separable, then $E$ is  J-reflexive.
\end{enumerate}
\end{corollary}
\Pr Point~\ref{it:c-b-c}. reformulates  Lemma~\ref{lem:sep-james}. For 
Point~\ref{it:sep}~:  if $E$ is separable, its dual ball, which is homeomorphic 
with a closed subset of $[-1,1]^{\IN}$,  is *weak sequentially compact~; 
so with Point~\ref{it:c-b-c},  $E$ has no extended triangular sequence. 
Besides the ``multiple'' version of  $\HB$ is provable in $\ZF$ 
for separable normed spaces, hence  $E$ has no triangular sequences either, whence 
it is J-reflexive.
\Rp

\medskip

\begin{remark} Using  Valdivia's theorem (see Question~\ref{ques:v-b-d-f-t}),
it follows from  $\DC$  that {\em Sup-reflexive spaces not containing isomorphic copies of
$\ell^1(\IN)$ do not have extended triangular sequences}.
\end{remark}

\medskip

\subsection{Spaces with *weak block compact dual ball} \label{subsec:*weak-db}
Say that a bounded sequence $(f_n)_{n \in \IN}$ of a normed space $E$
is {\em equivalent to the canonical basis of $\ell^1(\IN)$} if there exists some real number $M>0$ satisfying
$M \sum_{n \in \IN} |\lambda_n| \le \ln \sum_{n \in \IN} \lambda_n f_{n}\rn$ for every 
$(\lambda_n)_n \in \ell^1(\IN)$~: if in addition, $\forall n \in \IN \; \norme{f_n} \le 1$, 
say that $(f_n)_{n \in \IN}$ is {\em $M$-equivalent} 
to  the canonical basis of $\ell^1(\IN)$.
\begin{theorem*} [{\bf Rosenthal's $\ell^1$-theorem}] Given a 
set $X$ and a bounded sequence
$(f_n)_{n \in \IN}$ of $\ell^{\infty}(X)$, there exists a subsequence of $(f_n)_{n \in \IN}$
which pointwise converges, or there exists a subsequence which is equivalent to the canonical basis
of $\ell^1(\IN)$.
\end{theorem*}

Rosenthal's Theorem  
is a choiceless consequence (see Kechris, \cite[p.135-136]{Ke})
of the following choiceless result (see for example Avigad, \cite{Avi}, 1996)~: 
\begin{theorem*} [{\bf Cohen, Ehrenfeucht, Galvin (1967)}] Every open subset
of $[\IN]^{\IN}$ (the set of
infinite subsets of $\IN$ endowed with the product topology) is Ramsey.
\end{theorem*}

\medskip

\Pr[Proof of Theorem~\ref{lem:*wdual}]  Seeking a contradiction,
assume that some sup-reflexive normed space $E$ is not J-reflexive, though  $\Gamma_{E'}$
is *weak block  compact.
Non J-reflexivity of $E$ yields
 some $\vartheta$-triangular sequence $(a_n,f_n)_{n \in \IN}$ with $\vartheta >0$.
Using $\HB$, extend each $f_n$ to some $\tilde f_n \in \Gamma_{E'}$.
Then, with Lemma~\ref{lem:sep-james}, no infinite subsequence of  $(\tilde f_n)_{n \in \IN}$ 
pointwise converges, so,
using Rosenthal's $\ell^1$-theorem, there exists some infinite subsequence
$(\tilde f_n)_{n \in A}$ and some  $M>0$ such that the bounded sequence
$(\tilde f_n)_{n \in A}$ is $M$-equivalent to
the canonical basis of $\ell^1(\IN)$.  Now, by *weak
block compactness of $\Gamma_{E'}$,  $(\tilde f_n)_{n \in A}$  has a normalized
block sequence $(b_n)_{n \in \IN}$  *weakly converging  to $0$.
Using Simons' inequality, there exists some finite convex combination
$g:=\sum_{i \in F} \lambda_i b_i$ of $(b_n)_{n \in \IN}$ such that
$\ln g \rn =\sup_{\Gamma_E} g \le \frac{M}{2}$~;
but, since the block sequence $(b_n)_{n \in \IN}$ is normalized, it is also $M$-equivalent to the canonical basis
of $\ell^1(\IN)$, hence  $\ln g \rn \ge M \sum_{i \in F} |\lambda_i|=M$~: the contradiction~!
\Rp

\medskip

\section{Extension of a  theorem by Hagler and Johnson} \label{sec:ha-jo-zfc}
\begin{notation} [\cite{Ha-Jo}] If $(b_n)_{n \in \IN}$ is a normalized block sequence of a sequence 
$(x_n)_{n \in \IN}$ of a real vector space,  we write $(b_n)_n \prec (x_n)_n$. \\
Given a set $X$, for every  bounded sequence $(f_n)_{n \in \IN}$ 
of $\ell^{\infty}(X)$, 
and every subset $K$ of $X$, let 
$$\delta_K(f_n)_n := \sup_K \limsup f_n$$  
$$ \varepsilon_K(f_n)_n := \inf \{\delta_K(h_n)_n : \; (h_n)_n \prec (f_n)_n\}$$
\end{notation}

\medskip

\begin{remark} If for every $n \in \IN$, $f_n[K]=f_n[-K]$, 
then $\delta_K (f_n)_n = 0$ if and only if 
$(f_n)_n$ pointwise converges to $0$ on  $K$.
Observe that if for some $n_0 \in \IN$,  $(h_n)_{n \ge n_0}$ is a normalized block sequence of $(f_n)_n$, 
then $\delta_K(h_n) \le \delta_K(f_n)$ and  $\varepsilon_K(f_n)_n \le \varepsilon_K(h_n)$. 
\end{remark}

\medskip


Given a metric space $(X,d)$, for every $x \in X$ and every real number $r >0$
we denote by $B(x,r)$ the open ball $\{y \in X : \; d(x,y)  < r\}$.

\begin{lemma}[{\bf quantifier permuting}] \label{lem:limsup}  
Let $(K,d)$ be a precompact metric space, $\lambda \in \IR_+^*$, 
and $(f_n)_{n \in \IN}$ be a sequence of $\lambda$-Lipschitz real mappings on $K$. 
If  $\delta_K(f_n)_n \le 1$ then, for every  $\varepsilon \in \IR^*_+$, there exists $N \in \IN$ satisfying  
$\forall n \ge N  \; sup_{K} f_n \le 1+\varepsilon$. 
\end{lemma}
\Pr Let $\eta \in ]0,\varepsilon[$. Given some $x \in K$, there exists some 
finite subset $F_x$ of $\IN$ satisfying $\forall n \in \IN \backslash F_x \; f_n(x) < 1 + \eta$~; thus,
denoting by $\rho$ the positive number $\frac{\varepsilon - \eta}{\lambda}$, for every 
$n \in \IN \backslash F_x$, for every $y \in B(x,\rho)$, $f_n(y) < 1 + \varepsilon$.
Now the precompact set $K$ is contained in a 
finite union of the form $\bigcup_{1 \le i \le N} B(x_k,\rho)$. Let $F$ be the finite set 
$\bigcup_{1 \le k \le n} F_{x_k}$.
Then, for every $y \in K$, given $k \in \{1..n\}$ such that $d(x_k,y) < \rho$,
for every $i \in \IN \backslash F$, 
$f_i(y) < f_i(x_k) + \lambda \rho \le 1 + \eta + ( \varepsilon - \eta) = 1+\varepsilon$. 
\Rp


\medskip

\begin{lemma} [{\cite[proof of Theorem~1]{Ha-Jo}}] \label{lem:diag_ha-jo} 
Given a set $X$ and a
bounded sequence $(f_n)_{n \in \IN}$ of $\ell^{\infty}(X)$,
there exists a normalized block sequence $(b_n)_{n \in \IN}$ of $(f_n)_n$ such that
$\varepsilon_X(b_n)_n = \delta_X(b_n)_n$. 
\end{lemma}
\Pr Diagonalization. Choose some normalized block sequence $(h^0_n)_{n \in \IN}$ of $(f_n)_n$ such that 
$\delta_X(h^0_n)_n \le \varepsilon_X(f_n)_n + \frac{1}{2^0}$, and then, for every $i \in \IN$,
inductively choose some 
normalized block sequence $(h^{i+1}_n)_{n \in \IN}$ of $(h^i_n)_n$ such that
$\delta_X(h^{i+1}_n)_n \le \varepsilon_X(h^i_n)_n + \frac{1}{2^{i+1}}$.  
For every $n \in \IN$, let $b_n := h_n^n$~: then $(b_n)_n \prec (f_n)_n$~; moreover,
 given a normalized block sequence $(k_n)_{n \in \IN}$ of $(b_n)_n$, for every $i \in \IN$,
$\delta_X(b_n)_n  \le \delta_X(h^{i+1}_n)_n \le \varepsilon_X(h^i_n)_n +  \frac{1}{2^{i+1}} \le 
\delta_X(k_n)_n +  \frac{1}{2^{i+1}}$
whence $\delta_X(b_n)_n \le \varepsilon_X(b_n)_n $. 
\Rp

Lemma~\ref{lem:diag_ha-jo} is valid in $\ZF$~: no choice is needed here 
since block  can be built with rational coefficients.

\medskip

\begin{notation} For every integer $n \in \IN$,  we denote by $\mathcal S_n$ the set $\{0,1\}^n$  
of finite sequences of $\{0,1\}$   with length $n$~; let   $\mathcal S$  be the set of all 
finite sequences of $\{0,1\}$.
  Given an infinite subset $A$ of $\IN$, we denote by $i \mapsto i_A$ the increasing 
  mapping from $\IN$ onto $A$. 
\end{notation}

 Say that a family $(A_{\sigma})_{\sigma \in \mathcal S}$ of infinite subsets of $\IN$
  is a {\em tree} (of subsets of $\IN$) if  for every $\sigma \in \mathcal S$, 
$A_{\sigma \frown 0}$ and $A_{\sigma \frown 1}$ are disjoint subsets of $A_{\sigma}$.

\medskip

\Pr[Proof of Theorem~\ref{theo:ha-jo-zfc}] We essentially follow the proof of Hagler and Johnson, 
extending it with the help
of Lemma~\ref{lem:limsup}. Assuming  that $E$  is a normed space, and that 
$(g_n)_{n \in \IN}$ is a bounded sequence of $E'$ without any normalized block sequence 
pointwise converging to $0$, we are to show that $E$ contains an 
{\em asymptotically isometric} copy of $\ell^1(\IN)$.
Using Lemma~\ref{lem:diag_ha-jo}, the sequence $(g_n)_{n \in \IN}$ has a   
normalized block sequence $(f_n)_{n \in \IN}$
satisfying  $\varepsilon_{\Gamma_E}(f_n)_n = \delta_{\Gamma_E}(f_n)_n >0$. 
Dividing each $f_n$ by $\varepsilon_{\Gamma_E}(f_n)_n$, we
may assume that  $\varepsilon_{\Gamma_E}(f_n)_n = \delta_{\Gamma_E}(f_n)_n=1$.
Let $(u_n)_{n \in \IN}$ be a sequence of $]0,1/3[$ decreasing to $0$~; for every $n \in \IN$, let 
$\varepsilon_n := \frac{u_n}{2^n}$ and let $\delta_n := 1-\varepsilon_n$~;
thus  $(\varepsilon_n)_{n \in \IN}$
also decreases  to $0$.
Using $\DC$, we will build a tree $(A_{\sigma})_{\sigma \in \mathcal S}$ rooted at $A_{\varnothing}:=\IN$, 
and a sequence $(\omega_n)_{n \ge 1}$ of $\Gamma_E$
satisfying 
for every $n \ge 1$, $\sigma=(\alpha_1, \dots, \alpha_n) \in \mathcal S_n$, 
and $i \in A_{\sigma}$, 
 $$\ln {f_i}{\restriction span\{\omega_1, \dots, \omega_n\}} \rn \le 1 + \varepsilon_{n}
 \text{ and }
 \begin{cases} 
 f_i(\omega_{n}) \ge 1-3 u_{n} & \text{if } \alpha_n = 1 \\
 f_i(\omega_{n}) \le -1+3 u_{n} & \text{if } \alpha_n = 0 
\end{cases} 
 $$
 Then, with
 \begin{align*}
&P_n := \{f \in E' : f(\omega_n) \ge  1-3 u_{n} \text{ and }  
\ln f{\restriction span\{\omega_1, \dots,  \omega_n\}} \rn \le 1 + \varepsilon_n\} \\ 
&Q_n := \{f \in E' : f(\omega_n) \le  -1+3 u_{n} \text{ and }  
\ln f{\restriction span\{\omega_1, \dots, \omega_n\}} \rn \le 1 + \varepsilon_n\}
\end{align*}
it will follow that 
$(P_n, Q_n)_{n \ge 1}$ is {\em independent} (for every disjoint finite subsets $F, G$ of $\IN \backslash \{0\}$,
$\bigcap_{n \in F}P_n \cap \bigcap_{n \in G}Q_n$ is non-empty), whence the sequence
$(\omega_n)_{n \ge 1}$ of $\Gamma_E$ is asymptotically isometric  to the canonical basis
of $\ell^1(\IN)$~: indeed, given  real numbers $\lambda_1, \dots, \lambda_n$, letting 
$f \in \bigcap_{\{i : \lambda_i >0\}} P_i \cap \bigcap_{\{i : \lambda_i < 0\}} Q_i$,  
\begin{multline} \notag
{\ln f{\restriction span\{\omega_1, \dots, \omega_n\}} \rn}  \ln \sum_{1 \le i \le n} \lambda_i \omega_i \rn 
\geq f \big(  \sum_{1 \le i \le n} \lambda_i \omega_i \big)\\  
\ge \big(\sum_{\{i : \lambda_ i >0\}} \lambda_i (1-3u_i) + \sum_{\{i : \lambda_ i <0\}} \lambda_i (-1+3u_i)\big) 
\ge  \sum_{1 \le i \le n} |\lambda_i|(1-3u_i)
\end{multline}
whence
$\ln \sum_{1 \le i \le n} \lambda_i \omega_i \rn  
\ge  \frac{1}{1+\varepsilon_n} \sum_{1 \le i \le n} |\lambda_i|(1-3u_i) 
\ge \sum_{1 \le i \le n} |\lambda_i| \frac{1-3u_i}{1+\varepsilon_i}$,
with the sequence 
$(\frac{1-3 u_i}{1+\varepsilon_i})_{i \in \IN}$ of $]0,1[$ converging to $1$. \par
  {\em Building $\omega_{n+1}$ and $(A_{\sigma})_{\sigma \in \mathcal S_{n+1}}$
  from $(A_{\sigma})_{\sigma \in \mathcal S_n}$ and $(\omega_i)_{1 \le i \le n}$.} 
  Given, for every $\sigma \in \mathcal S_n$,  two infinite disjoint subsets 
 $L_{\sigma}$ and $R_{\sigma}$ of $A_{\sigma}$, 
  consider the normalized block sequence $(h^n_i)_{i \in \IN}$ with 
  $h^n_i := \frac{1}{2^n} \sum_{\sigma \in \mathcal S_n} 
  \frac{f_{i_{R_{\sigma}}} - f_{i_{L_{\sigma}}}}{2} $.  
  Since $\delta_{\Gamma_E}(h^n_i)_i \ge 1$,   
  there is some $\omega_{n+1} \in \Gamma_E$ satisfying
  $\limsup_i h^n_i(\omega_{n+1}) > \delta_{n+1}$ and in particular,
  the set   $J:=\{i \in \IN \; : h^n_i(\omega_{n+1}) > \delta_{n+1}\}$ is infinite. 
  Since  the closed unit ball 
  $K$ of the finite dimensional space $span\{\omega_1, \dots, \omega_{n+1}\}$ is compact and 
  $\delta_K(f_i)_i \le \delta_{\Gamma_E}(f_i)_i \le 1$, 
    Lemma~\ref{lem:limsup} implies the existence of some $N_0 \in \IN$ satisfying 
  $\forall i \ge N_0 \; 
  \ln {f_i}{\restriction span\{\omega_1, \dots, \omega_{n+1}\}} \rn \le 1+\varepsilon_{n+1}$.
  Let   $J'  := \{ i \in J :  \; i \ge N_0\}$. 
  Now, given any  $\sigma \in \mathcal S_n$,  for every
  $i \in J'$, notice that, below, $i_{R_{\tau}}, i_{L_{\tau}} \ge i \ge N_0$~:
    \begin{equation} \notag 
  \begin{split}
  \frac{f_{i_{R_{\sigma}}}(\omega_{n+1}) - f_{i_{L_{\sigma}}}(\omega_{n+1})}{2} 
  &=     2^{n} h^n_i(\omega_{n+1}) - \sum_{\tau \in \mathcal S_n, \tau \neq \sigma} 
  \frac{f_{i_{R_{\tau}}}(\omega_{n+1}) - f_{i_{L_{\tau}}}(\omega_{n+1})}{2} \\
  &\ge     2^{n} \delta_{n+1} - (2^n-1) (1+\varepsilon_{n+1}) \\
  &= 2^{n} (1-\varepsilon_{n+1}) - (2^n-1) (1+\varepsilon_{n+1}) \\
  &=1  - (2^{n+1}-1) \varepsilon_{n+1}  \ge    1  - 2^{n+1} \varepsilon_{n+1} =1-u_{n+1} 
   \end{split}
   \end{equation}
   whence 
   \begin{equation} \notag
  \begin{split}
   f_{i_{R_{\sigma}}}(\omega_{n+1}) &\ge 2(1-u_{n+1}) + f_{i_{L_{\sigma}}}(\omega_{n+1}) \\
   &\ge 2(1-u_{n+1}) -(1+\varepsilon_{n+1}) \\
   &=1-u_{n+1}(2+\frac{1}{2^{n+1}}) \ge 1-3u_{n+1}
   \end{split}
   \end{equation}
    \noindent likewise,   $f_{i_{L_{\sigma}}}(\omega_{n+1}) \le - 1+3 u_{n+1}$~;  
	then, let 
  $A_{\sigma \frown 0} := \{i_{L_{\sigma}} : i \in J'\}$, and 
   $A_{\sigma \frown 1} := \{i_{R_{\sigma}} : i \in J'\}$.
   \Rp


\medskip

\section{No asymptotically isometric copy of $\ell^1(\IN)$ in sup-reflexive spaces} 
\label{sec:noasyiso}
Notice that any asymptotically isometric  copy $(a_n)_{n \in \IN}$ of
$\ell^1(\IN)$ in a normed space $E$  is linearly independent 
(indeed, it is a Schauder sequence of $E$).

\medskip
\Pr[Proof of Theorem~\ref{theo:asy-l1}] 
Assume the existence of some sequence $(a_n)_{n \in \IN}$ of  $\Gamma_E$, 
asymptotically isometric with the canonical
basis of $\ell^1(\IN)$, witnessed by a sequence of  coefficients $(\delta_i)_{i \in \IN}$ of $]0,1]$ 
converging to $1$.
Let $V:= {span}\{a_n : n \in \IN\}$. 
 For every $n \in \IN$, consider  the linear
 mapping  $g_n : V \to \IR$ such that $g_n(a_i)=-\delta_i $ if $i <n$ and
 $g_n(a_i)=\delta_i$ if $n \le i$~; then, for every $(\lambda_i)_{i \in \IN} \in \IR^{(\IN)}$
 ($\lambda_i=0$ for all but finitely many $i$'s),
    $|g_n(\sum_{j=0}^{\infty} \lambda_j a_j )| =
 |\sum_{j <n} -\lambda_j \delta_j + \sum_{j \ge n} \lambda_j \delta_j|
 \le \sum_{j=0}^{\infty} |\lambda_j| \delta_j \le 
 \ln \sum_{j=0}^{\infty} \lambda_j a_j \rn $
 whence $g_n$ is continuous and $\ln g_n \rn \le 1$~; besides
 for every integer $i \ge n$, $g_n(a_i)=\delta_i$,
 whence $\lim_{i \to +\infty}g_n(a_i)=1$~; so $\ln g_n \rn=1$.
  Using $\HB$,  for each $n \in \IN$, extend $g_n$ to some 
 $\tilde g_n \in S_{E'}$. Let $W$ be the vector subspace of elements
 $x \in E$ such that $\big(\tilde g_n(x)\big)_{n \in \IN}$ converges. Then 
the linear mapping  $g:=\lim_n \tilde g_n$ is continuous with norm $\le 1$ on $W$~;
extend it to some element $\tilde g \in \Gamma_{E'}$. 
Now consider some sequence 
$(\alpha_i)_{i \in \IN}$ of $]0,1[$ such that $\sum_{i \in \IN} \alpha_i=1$ and
let $h:=\sum_{k \in \IN} \alpha_k \tilde g_k - \tilde g$. Clearly, $\ln h \rn \le 2$. 
Besides, for every $n \in \IN$,
$h(a_n)=\sum_{k \le n } \alpha_k \delta_n - \sum_{k > n } \alpha_k \delta_n +  \delta_n=
2 \delta_n \sum_{k \le  n } \alpha_k$,
thus   $\lim_n h(a_n)=2$.
So $\ln h \rn =2$. By sup-reflexivity of $E$,
let $u \in \Gamma_E$ such that $h(u)=2$. Observe that $\tilde g(u)=-1$,
and  for every $k \in \IN$,
$\tilde g_k(u)=1$ (notice that  for each $k$, $\alpha_k \neq 0$)~;
now $u \in W$ and $g(u)=\lim_k \tilde g_k(u)=1$, contradicting $\tilde g(u)=-1$~!
\Rp

\medskip

\section{Comments in set-theory  without choice $\ZF$} \label{sec:AC}
Recall the Axiom of Choice~:
\begin{axiom}[$\AC$] {\em Given a family $(A_i)_{i \in I}$ of non-empty sets, there exists a mapping
$f : I \to \cup_{i \in I} A_i$ satisfying $f(i) \in A_i$ for every $i \in I$.}
\end{axiom}

\noindent and  the two following weak forms  of the Axiom of Choice~:
\begin{axiom}[$\DC$] {\bf (Dependent Choices)} 
{\em Given a non-empty set $E$ and a binary relation $R \subseteq E \times E$ satisfying
$\forall x \in E \; \exists y \in E \; x R y$, there exists a sequence $(x_n)_{n \in \IN}$ of $E$
satisfying $\forall n \in \IN x_n R x_{n+1}$.}
\end{axiom}

\begin{axiom}[$\HB$] {\bf (Hahn-Banach)} 
{\em Given a real vector space $E$, a sublinear mapping $p : E \to \IR$, 
a vector subspace $F$ of $E$
and a linear mapping $f : F \to \IR$ satisfying $f \le p{\restriction F}$, 
there exists a linear mapping
$\tilde f : E \to \IR$ extending $f$ and satisfying $\tilde f \le p$.}
\end{axiom}

It is known (see~\cite{Ho-Ru}) that $\ZF \mathbf + \DC$ does not prove $\HB$,  that 
$\ZF \mathbf + \HB$ does not prove $\DC$, and that $\ZF \mathbf + \HB \mathbf + \DC$ does not prove 
$\AC$. The axiom $\HB$ is known to be equivalent to its {\em multiple form} (see~\cite{Ho-Ru})~: 
{\em Given a family $(E_i)_{i \in I}$ of real vector spaces,
a family $(p_i)_{i \in I}$ of sublinear mappings $p_i : E_i \to \IR$, 
a family $(F_i)_{i \in I}$ of vector subspaces $F_i \subseteq E_i$
and a family $(f_i)_{i \in I}$ of linear mappings $f_i : F_i \to \IR$ satisfying 
$f_i \le p_i{\restriction F_i}$ for every $i \in I$,
there exists a family $(\tilde f_i)_{i \in I}$ of linear mappings
$\tilde f_i : E_i \to \IR$ extending $f_i$ and satisfying $\tilde f_i \le p_i$
for every $i \in I$.}

\medskip  

\begin{question} All notions of reflexivity we reviewed  till now 
(J-reflexivity, onto-reflexivity, convex-reflexivity,
compact-reflexivity, Smulian-reflexivity, sup-reflexivity, \dots -see~\cite{Mo04}-) 
are equivalent in $\ZF \mathbf + \DC \mathbf + \HB$.
Is there some ``classical'' notion of reflexivity which is not equivalent to ``sup-reflexivity'' in
$\ZF \mathbf + \DC \mathbf + \HB$~?
\end{question}

\begin{question} In $\ZF$ (see Corollary~\ref{cor:sep} of Section~\ref{sec:dual-ball-si}),
separable sup-reflexive Banach spaces are J-reflexive. More generally,
all notions of reflexivity for {\em separable} Banach spaces
we reviewed  till now are equivalent in $\ZF$ either to onto-reflexivity or to J-reflexivity, 
the former being weaker than the latter, and 
these two notions being not equivalent since there exists a model of $\ZF \mathbf + \DC$ 
(see~\cite{Pi-So}) where the
continuous dual of $\ell^{\infty}(\IN)$ is $\ell^1(\IN)$~: in such a model, the
(non J-reflexive, separable) Banach space $\ell^1(\IN)$  is onto-reflexive.
Is there some ``classical'' notion of reflexivity which, for separable Banach spaces,
is  equivalent in $\ZF$, neither to J-reflexivity 
nor to onto-reflexivity~?
\end{question}

\medskip

\bibliographystyle{abbrv}
\bibliography{envoi.bbl}

\begin{thebibliography}{10}

\bibitem{Avi}
Avigad.
\newblock {A new proof that open sets are Ramsey}.
\newblock {\em Arch. Math. Logic}, 37:235--240, 1998.

\bibitem{Bou-Die}
Bourgain and Diestel.
\newblock Limited operators and strict cosingularity.
\newblock {\em Math. Nachrichten}, 119:55--58, 1984.

\bibitem{Bou-Fre-Tal}
{Bourgain, Fremlin, Talagrand}.
\newblock {Pointwise compact sets of Baire-measurable functions}.
\newblock {\em {Amer. J. Math.}}, 100:845--886, 1978.

\bibitem{DGZ}
R.~Deville, G.~Godefroy, and V.~Zizler.
\newblock {\em {Smoothness and renormings in Banach spaces.}}, volume~64.
\newblock New York: John Wiley \& Sons, Inc., 1993.

\bibitem{Ha-Jo}
{Hagler, J. and Johnson, W.B.}
\newblock {On Banach spaces whose dual balls are not weak* sequentially
  compact}.
\newblock {\em Isr. J. Math.}, 28(4):325--330, 1977.

\bibitem{Ha-Od}
{Hagler, J. and Odell, J.}
\newblock {A Banach space not containing $\ell^1$ whose dual ball is not weak*
  sequentially compact}.
\newblock {\em {Illinois Journal of Mathematics}}, 22(2):290--294, 1978.

\bibitem{Ho-Ru}
P.~Howard and J.~E. Rubin.
\newblock {\em {Consequences of the Axiom of Choice.}}, volume~59.
\newblock American Mathematical Society, Providence, RI, 1998.

\bibitem{Ja64c}
{James, R.C.}
\newblock {Characterizations of reflexivity}.
\newblock {\em Stud. Math.}, 23:205--216, 1964.

\bibitem{Ja64b}
{James, R.C.}
\newblock {Weak compactness and reflexivity}.
\newblock {\em Isr. J. Math.}, 2:101--119, 1964.

\bibitem{Ja72}
{James, R.C.}
\newblock Reflexivity and the sup of linear functionals.
\newblock {\em Isr. J. Math.}, 13:289--300, 1972.

\bibitem{Ke}
A.~S. Kechris.
\newblock {\em {Classical descriptive set theory.}}
\newblock Springer-Verlag, Berlin, {GTM 156} edition, 1994.

\bibitem{Mo04}
M.~Morillon.
\newblock {James sequences and Dependent Choices}.
\newblock {\em Math. Log. Quart.}, 51(2), 2005.

\bibitem{oj98b}
E.~Oja.
\newblock {A proof of the Simons inequality}.
\newblock {\em Acta et Commentationes Universitatis Tartuensis de Mathematica},
  2:27--28, 1998.

\bibitem{Oj98}
E.~Oja.
\newblock {A short proof of a characterization of reflexivity of James.}
\newblock {\em Proc. Am. Math. Soc.}, 126(8):2507--2508, 1998.

\bibitem{Pi-So}
D.~Pincus and R.~M. Solovay.
\newblock {Definability of measures and ultrafilters.}
\newblock {\em J. Symb. Log.}, 42:179--190, 1977.

\bibitem{Pr}
Pryce.
\newblock Weak compactness in locally convex spaces.
\newblock {\em {Proc. Amer. Math. Soc.}}, 17:148--155, 1966.

\bibitem{Si72}
S.~Simons.
\newblock A convergence theorem with boundary.
\newblock {\em {Pacific J. Math.}}, 40:703--721, 1972.

\bibitem{Va}
Valdivia.
\newblock Fr\'echet spaces whith no subspaces isomorphic to $\ell^1$.
\newblock {\em Math. Japon.}, 38:397--411, 1993.

\end{thebibliography}
\end{document}